\documentclass[11pt,a4paper]{article}

\usepackage{amsmath}
\usepackage{amssymb}
\usepackage{verbatim}
\usepackage{amsthm}
\usepackage{mathrsfs}
\usepackage{multirow}
\usepackage{stackengine}
\usepackage{scalerel}
\usepackage{graphicx}
\usepackage{mathtools}
\usepackage[colorlinks,pdfdisplaydoctitle,linkcolor=blue,citecolor=blue]{hyperref}   

\usepackage{caption}
\usepackage{subcaption}
\usepackage{url}
\usepackage{epstopdf}
\usepackage{pgf,tikz}
\usepackage{bbm}
\usepackage{extarrows}
\usetikzlibrary{arrows}
\usepackage[ruled,lined]{algorithm2e}
\usepackage{titling}

\usepackage{dsfont}
\usepackage{listings}
\usepackage{bbm}
\usepackage{color}
\usepackage{dsfont}
\usepackage{enumerate}
\DeclareMathOperator*{\argmin}{arg\,min}

\usepackage{tabularx}

\renewcommand{\Delta}{\triangle}

\definecolor{darkblue}{rgb}{0,0,0.7}

\usepackage{soul}

\definecolor{darkgreen}{rgb}{0.01,0.75,0.24}
\newcommand{\bbR}{\mathbb{R}}

\def \Ee[#1]{\mathcal{E}^{\text{{#1}}}}
\def\R{\mathbf{R}}

\def\pa[#1,#2]{\frac{\partial {#1}}{\partial {#2}} }

\def\idom[#1,#2,#3]{\int_{#1}\hspace{1pt} {#2} \hspace{1pt} \text{d}{#3}}
\def\res[#1,#2]{\left.{#1}\right|_{#2}}

\def\var[#1,#2]{\langle \delta \mathcal{E}^{\text{{#1}}}({#2}),v\rangle}
\def\vars[#1,#2,#3]{\langle \delta^2\mathcal{E}^{\text{{#1}}}({#2})v,{#3}\rangle}
\def\vard[#1,#2,#3,#4]{\langle \delta\mathcal{E}^{\text{{#1}}}({#2})-\delta\mathcal{E}^{\text{{#3}}}({#4}),v\rangle}

\newcommand{\Cov}[1]{\mathrm{Cov}\left[#1\right]}

\newcommand{\C}{\mathcal{C}}

\newcommand{\be}{\begin{equation}}
\newcommand{\en}{\end{equation}}
\newcommand{\ben}{\begin{equation*}}
\newcommand{\enn}{\end{equation*}}
\newcommand{\bea}{\begin{aligned}}
\newcommand{\ena}{\end{aligned}}

\def\ba#1\ena{\begin{align}#1\end{align}}
\def\ban#1\enan{\begin{align*}#1\end{align*}}

\usepackage[colorinlistoftodos,bordercolor=orange,backgroundcolor=red!25,linecolor=orange,textsize=scriptsize]{todonotes}

\newcommand{\houcine}[1]{\todo[inline]{\textbf{Houcine: }#1}}

\newcommand{\neil}[1]{\todo[inline]{\textbf{Neil: }#1}}

\theoremstyle{plain}
\newtheorem{thm}{Theorem}[section]
\newtheorem{theorem}{Theorem}[section]

\newtheorem{Lem}[thm]{Lemma}

\newtheorem{lemma}[thm]{Lemma}
\newtheorem{cor}[thm]{Corollary}

\newtheorem{ass}{Assumption}

\newtheorem{remark}[thm]{Remark}
\newtheorem{rem}[thm]{Remark}

\numberwithin{equation}{section}

\newcommand{\X}{\mathcal{X}}
\newcommand{\Y}{\mathcal Y}

\newcommand{\op}{F}
\newcommand{\norm}[2]{\left\Vert #1\right\Vert_{#2}}
\newcommand{\ip}[2]{\left\langle #1, #2\right\rangle}
\newcommand{\taylor}[2]{\op\left(#1\right) + \op'\left[#1\right]\left(#2-#1\right)}

\newcommand{\udag}{u^\dagger}
\newcommand{\ydag}{y^\dagger}
\newcommand{\cS}[1]{\mathcal S_n \left(#1\right)}
\newcommand{\cSt}[1]{\tilde{\mathcal S}_n \left(#1\right)}
\newcommand{\cT}[1]{\frac{1}{2}\|#1- \ydag\|^2_{\mathcal{Y}}}
\newcommand{\err}{\mathrm{err}_n}
\newcommand{\tildeerr}{\widetilde{\mathrm{err}}_n}

\newcommand{\Ctc}{C_{\mathrm{tc}}}

 \def\emailname{E-mail}%
\def\email#1{\emailname: #1}

\begin{document}

\title{A Stochastic  Iteratively Regularized Gauss-Newton  Method}

\author{
El Houcine Bergou\thanks{Mohammed VI Polytechnic University, Ben Guerir, Morocco. 
\email{{\tt elhoucine.bergou@um6p.ma}}}
\and
Neil K. Chada\thanks{Department of Mathematics, City University of Hong Kong, 83 Tat Chee Avenue, Kowloon Tong, Hong Kong. \email{{\tt neilchada123@gmail.com}}
  }      
  \and
Youssef Diouane\thanks{GERAD and Department of Mathematics and Industrial Engineering, Polytechnique Montréal, Montréal, Canada. \email{{\tt youssef.diouane@polymtl.ca}}
  }  
}

\date{\today}
\maketitle

\begin{abstract}
This work focuses on developing and motivating a stochastic version of a well-known inverse problem methodology. Specifically, we consider the iteratively regularized Gauss-Newton method, originally proposed by Bakushinskii for infinite-dimensional problems. Recent work  have extended this method to handle sequential observations, rather than a single instance of the data, demonstrating notable improvements in reconstruction accuracy. In this paper, we further extend these methods to a stochastic framework through mini-batching, introducing a new algorithm, the stochastic iteratively regularized Gauss-Newton method (SIRGNM). Our algorithm is designed through the use randomized sketching. We provide an analysis for the SIRGNM, which includes a preliminary error decomposition and a convergence analysis, related to the residuals. We provide numerical experiments on a 2D elliptic PDE example. This illustrates the effectiveness of the SIRGNM, through maintaining a similar level of accuracy while reducing on the computational time.
 \\\\
 \textbf{Subject Class:} 65N21, 65C35, 65K10, 93E24.
\end{abstract}




\section{Introduction}

The field of optimization \cite{BV04,YN04,NW06} is of crucial importance many areas of applied mathematics, which is concerned with minimizing functions, or functionals of interest. Popular examples include control theory, calculus of variations, numerical analysis and others. This is particularly apparent with its recent success in the emerging fields of machine learning, data science and deep learning \cite{BCN18,SCZ20,WWJ16}. However a particular application, which is of interest for this work, is parameter estimation associated with differential equations. This mathematical discipline is commonly known as inverse problems \cite{LP13,AMS10,AT87}, which is concerned with the recovery of some quantity, or parameter, $u^{\dagger} \in \mathcal{X}$ of interest from some noisy observations $y^{\dagger} \in \mathcal{Y}$, i.e.
\begin{equation}
\label{eq:ip}
y^{\dagger} = F(u^{\dagger}),
\end{equation}
with $F:\mathcal{X} \rightarrow \mathcal{Y}$ being the forward operator associated with the targeted problem.
In practice, due to the unavoidable presence of observational noise in real applications,
equation~\eqref{eq:ip} must be replaced by
\begin{equation}
\label{eq:ip:noise}
y^{\delta} = F(u^{\dagger}) + \eta, \quad \eta \sim \mathcal{N}(0,\delta^2I),
\end{equation}
 where $y^{\delta} \in \mathcal{Y}$ represents the corrupt noisy data, $\eta$ takes the form of additive Gaussian nose, and $\delta>0$. Traditionally inverse problems are solved using variational techniques, motivated from regularized least-squares, where the solution of \eqref{eq:ip:noise} is estimated by the following minimizer
\begin{equation} \label{eq:ustar}
u^*\in \argmin_{u \in \mathcal{X}} \frac{1}{2} \Big\| y^{\delta} - F(u)\big\|^2_{\mathcal{Y}} + \alpha \|u\|^2_{\mathcal{X}},
\end{equation}

where $\alpha$ is a step-size rule, and $\|u\|^2_{\mathcal{X}}$ acts as a penalty term. This is a well-used, and popular choice of penalty known as Tikhonov regularization. There has been huge developments in the area of developing 
regularized least-squares solvers for inverse problems, which are ill-posed, known as iterative regularization methods~\cite{EHN96,KNS08}.
Well-known algorithms based on this include the regularized Levenberg-Marquardt method (LMM), $\nu$-methods and the Landweber iteration \cite{BGV16}. Our focus in this work will be another method, which has connections to the above stated methods which is the iterated regularized Gauss--Newton method (IRGNM) proposed by Baksushinkii~\cite{ABB92,BK04}, which aims to mimic the \textit{classical} iterated Gauss--Newton method, but for ill-posed problems. Since its development it has seen significant advances both related to theory and applications. Such example are the development of error bounds in different function space settings, demonstrating convergence and applications to data assimilation-based methods, such as when one has sequential data~\cite{BHM09,BNO97,CIL22,HW13,QJ00,QJ11,FW15}.  Despite these important directions and works, there has been somewhat a 
discontinuity between classical optimization understanding, and development, of algorithms, and those related to
inverse problems. In particular a recent attraction has been that of stochastic gradient methods \cite{BCN18,BV04,RM51}, which have proven popular due to their cost-effectiveness, compared to full-gradient methods. This is related to not using the full derivative information, but rather a random subset of it. The most well-known example of this includes the stochastic gradient descent. In the context of least-squares optimization this has also been done, for the non-regularized Gauss--Newton method, the LMM and evolutionary algorithms. However this has not been very well-understood in the context of iterative regularization methods for inverse problems. Recent work has looked at understanding stochastic iterative methods for inverse problems, however much of the analysis is assumed to be in the linear and finite-dimensional setting \cite{HBD21,JC24,JL24}. Therefore this acts as our primary motivation to apply these techniques to iterative regularization methods, such as the IRGNM. 

Specifically our form of stochasticty we apply is based on random projection methods which have been introduced in \cite{GR15,RM16,RT20}, which is also known as sketching. These forms of projection has been shown to be easy to apply with desirable convergence properties. It is also well-known that sketching, in this form, is equivalent to using mini batching within optimization. Another motivation from the development of a stochastic version, is the derivation of convergence analysis which has not been considered before. In the classical IRGNM this has not been considered. We will consider both the IRGNM and a recent dynamic version (with multiple observations \cite{CIL22}), which are both presented in the infinite-dimensional setting. 
\subsection{Contributions}
Our highlighting contributions from this paper are provided below:
\begin{itemize}
    \item We propose a stochastic version of the IRGNM,
    which multiplies a projection operator by the forward operator. This technique is related
    to the sketching. Our method is entitled the stochastic iterated regularizaed Gauss-Newton method (SIRGNM). The addition of stochastic gradients is considered for both the classic and sequential/dynamic IRGNM.
    \item We present an analysis of the convergence properties of SIRGNM, which extends the well-established properties of the deterministic IRGNM. Our analysis is based on a projection operator for the sketching, where we present an error decomposition and a convergence analysis, which includes the residuals.
    \item 
    We evaluate the efficiency of our novel  approach through experimentation on a 2D partial differential equation (PDE)-based inverse problem. Our PDE of choice is referred to as Darcy, or groundwater flow, concerned with modelling fluid flow in a porous medium underground. We demonstrate the efficiency of exploiting a stochastic version while attaining a similar order of accuracy. 
\end{itemize}

\subsection{Outline}
The outline of this work is as follows. In Section \ref{sec:back} 
we provide a primer on the necessary material required for the manuscript.
This includes a review on the IRGNM, while discussing the aspects of projection, and sketching. This will lead onto Section \ref{sec:sirgnm} where we present
our new developed methodology of the stochastic method SIRGNM. We aim to carry out an analysis related to convergence, which will be presented in Section \ref{sec:theory_det}. Numerical results are presented in Section \ref{sec:num} demonstrating the effectiveness of the our proposed method. Finally we present an overview of our work, while mentioning fruitful areas of future work, to be conducted in Section \ref{sec:conc}.

\section{Background material}\label{sec:back}

In this section we briefly recall our motivating algorithm, which is the iterative regularized Gauss--Newton method and its dynamic counterpart, i.e. the dIRGNM.

Without loss of generality and to avoid unnecessarily notation we assume $y^\delta = 0$ in \eqref{eq:ustar}. 
Hence, our target in this paper reduces to solving the following least squares problem:
\begin{equation} \label{prob:optim}
  \min_{u\in \mathcal{X}} f(u) := \frac{1}{2} \| F(u)\|_{\mathcal{Y}}^2 + \alpha \|u\|^2_{\mathcal{X}},
\end{equation}
A popular and well-known methodology for solving \eqref{prob:optim} is the iterated regularizaed Gauss--Newton method (iRGNM). The variational form based on the IRGNM is defined as:
\begin{equation}
\label{eq:var_gn}
u_{n+1} := \argmin_{u \in \mathcal{X}}\left[\norm{F\left(u_n\right) + F'\left[ u_n\right] \left(u- u_n\right)}{\Y}^2  + \alpha_n \norm{u- u_0}{\X}^2\right],
\end{equation}
where $\alpha_n>0$ is a the step-size (which os iteration dependent).  The control of the parameter $\alpha_n$ of the iterations is often essential for the convergence of the iterative process. 
In our setting, $u_0 \in \X$ represents the apriori information on $u_n$.  
If $u_0  =u_n$, it corresponds to the classical Levenberg-Marquardt (LM) method. $u_0$ can be also set to some initial guess that does not depend on the iterative process~\cite{CCS21}. $F'\left[u\right]$ is the  Fr\'{e}chet  derivative of $F$ at $u$. Without loss of generality,  a covariance operator can be included to  scale the term $\norm{u- u_0}{\X}^2$. Namely, the latter term can be replaced with $\norm{\C^{-1/2}(u- u_0)}{\X}^2$ where $\C$ is a covariance operator, which we will discuss in more details in the numerical section.

In the work of \cite{CCS21}, the authors designed an alternative version for \eqref{eq:var_gn}, which takes into consideration random sequential data, motivated from the field of data assimilation. 
Specifically, the observational model with \textit{sequential noisy observations} has the form
\begin{equation}
\label{eq:model_dyn}
Y_i = \op(\udag) + \sigma \xi_i,\quad  i=1,2,\ldots,
\end{equation}
where now the averaged observations, $Z_n$, of \eqref{eq:model_dyn} is defined as
\begin{align}
\label{eq:avg_noise}
Z_n = n^{-1} \sum_{i=1}^nY_i = \op\left(\udag\right)+ \frac{\sigma}{n} \sum_{i=1}^n \xi_i.
\end{align}
The motivation behind considering \eqref{eq:avg_noise} is that the covariance of the data decreases as $i$ increasing, unlike the fixed data choice of $Y$ which is $\sigma$, i.e.
$$
\textrm{Cov}{Z_n} = \Cov{\frac{\sigma}{n} \sum_{i=1}^n \xi_i} = \frac{\sigma^2}{n^2} \sum_{i=1}^n \textrm{Cov}{\xi_i} = \frac{\sigma^2}{n} \text{id}{\mathcal{Y}}.
$$
For the negative log-likelihood functional of the normal distribution, it replaces the term
$g:= {F\left(u_n\right) + F'\left[ u_n\right] \left(u- u_n\right)}$ by
\begin{align}\label{eq:S1}
\mathcal{S}_n(F\left(u_n\right) + F'\left[ u_n\right] \left(u- u_n\right)) &:= \frac12 \norm{F\left(u_n\right) + F'\left[ u_n\right] \left(u- u_n\right)}{\Y}^2 \\
& - \ip{F\left(u_n\right) + F'\left[ u_n\right] \left(u- u_n\right)}{Z_n}, \nonumber
\end{align}
where $Z_n \in \mathcal{Y}$ is observation dependent term, typically $Z_n$ is set as the average of a selected subset of sequential observations.

%
Hence, instead of \eqref{eq:var_gn}, one uses the following method modification:
\begin{align}
u_{n+1} := \argmin_{u \in \X} \bigg[\mathcal S_n \left(\taylor{u_n}{u}\right) + \alpha_n \norm{u - u_0}{\X}^2\bigg]. \label{eq_IRGNM}
\end{align}

\begin{remark}
\label{rem:1}
An important question is why do we consider the negative log-likelihood defined in \eqref{eq:S1} compared to 
$\|g\|_{\Y}^2$? The reason for this is because, in infinite dimensions, the quantity $\|g\|_{\Y}^2$ is almost-surely infinite, with respect to the data $Y$. This is because Gaussian measures in infinite dimensions do not lie in the associated Cameron-Martin space $\mathrm{Im}(\mathcal{Y}^{1/2})$.
This is also the case because $Y$ as a function rather than a finite-dimensional vector. Therefore what we do is subtract off the ``infinite" part in \eqref{eq:S1}. 
\end{remark}


Alternatively we can express the minimization procedure of \eqref{eq_IRGNM} in terms of the first order optimality condition as
\begin{equation}
\label{eq:gn}
u_{n+1} =  u_n + (F'\left[u_n\right]^*F'\left[ u_n\right]+ \alpha_n \mbox{Id}_{\X})^{-1} \left(  F'\left[u_n\right]^*\left(Z_n -F\left( u_n\right)\right)  + \alpha_n \left(u_0- u_n\right)\right),
\end{equation}
where $F'\left[u_n\right]^* : \Y \to \X$ is the adjoint operator of $F'\left[u_n\right] : \X \to \Y$. For computational efficiency, we can then use
Woodbury lemma for \eqref{eq:gn} yielding
\begin{align}\label{eq:gn:2}
u_{n+1}=u_{0}+\op'[u_{n}]^*(\op'[ u_{n}]\op'[ u_{n}]^*+\alpha_{n}\mbox{Id}_{\X})^{-1}\Big(Z_n-\op(u_{n})-\op'[u_{n}](u_0- u_{n})\Big).
\end{align}
An overview of the iterated regularized Gauss-Newton method (\texttt{IRGNM}) is provided in Algorithm \ref{alg:cIRGNM}. We note that we now take the notation for the maximum number of iterations to $M$.

\begin{algorithm}[h]
\caption{Iterated regularized Gauss-Newton method}
\label{alg:cIRGNM}
\SetKwInOut{Input}{inputs}
 \SetKwInOut{Output}{output}
    \Input{$u_0 \in \X$, $\alpha_0>0$ (an initial value for the regularisation), and $M$ (a maximum number of iterations).}


\For{$n=0,\ldots,M-1$}{
Select  $Z_n \in \Y$.\\
Set the ${\mathcal{S}}_n$ as defined in (\ref{eq:S1}).\\
Compute $u_{n+1}$, i.e., $u_{n+1} := \displaystyle \argmin_{u \in \X} \Bigg[{\mathcal S}_n \left(\taylor{u_{n}}{u}\right) + \alpha_n \norm{u - u_0}{\X}^2\Bigg].$\\
Update the regularization parameter $\alpha_{n+1}$.
}
    \Output{${u}_{M}$.}
\end{algorithm}
\section{Stochastic IRGNM}
\label{sec:sirgnm}
The recent success and popularity of variational, or optimization schemes, has been evident in the areas of machine learning and data science, where one aims to optimize a function or objective function. In particular this has been the case primarily due to stochastic optimization methods, which are known to have faster convergence properties while taking less computational time than evaluating the full gradient. The disadvantage of these methodologies is generally the accuracy is slightly worse. Given this, our aim is to incorporate ideas from stochastic gradient descent applied to the IRGNM, in the context of solving \eqref{eq:ip}.

\subsection{Towards a general stochastic framework for IRGNM}
To initiate this we can take two particular approaches in doing so. The first is based on a sketching operator, that projects the problem on a \textit{lower-dimensional} space, i.e., one uses
 $\mathcal{P}\left[ F(u)\right]$ as a stochastic estimation of $F(u)$, where $\mathcal{P}$ is the stochastic sketch operator that we assume is unbiased 


As a result we can include the modified operator $\mathcal{P}_n$ (depending on the iteration index $n$) in the variational problem as
\begin{align}
u_{n+1} := \argmin_{u \in \X} \bigg[\tilde{\mathcal S}_n \left(\taylor{u_n}{u}\right) + \alpha_n \norm{(u - u_0)}{\X}^2\bigg]. \label{eq:var_sgn}
\end{align}

The operator $\tilde{\mathcal S}_n$ is defined as follows: 
\begin{equation}
\label{eq:S}
\tilde{\mathcal S}_n \left(g\right)  
:= \frac{1}{2}\|\mathcal{P}_n\left[g\right]\|^2_{\mathcal{Y}_p}- \left< Z_n,  \mathcal{P}_n\left[g\right]\right>.
\end{equation}
where $\mathcal{Y}_p$ is set as the projected observations space. The vector point $Z_n \in \mathcal{Y}_p$ is selected by the user.
The proposed framework is detailed in Algorithm~\ref{alg:SIRGNM}. We note that the update of the regularization parameter is purposefully left unspecified at this stage. From now on, we will refer to our proposed approsch as stochastic iterated regularized Gauss-Newton method (\texttt{SIRGNM}).

\begin{algorithm}[h!]
\caption{Stochastic iterated regularized Gauss-Newton method}
\label{alg:SIRGNM}
\SetKwInOut{Input}{Inputs}
 \SetKwInOut{Output}{Output}
     \Input{$u_0 \in \X$, $\alpha_0>0$ (an initial value for the regularisation), and $M$ (a maximum number of iterations).}


\For{$n=0,\ldots,M-1$}{
Select  $\mathcal{P}_n$ a stochastic sketching operator.\\
Select  $Z_n \in \Y_p$.\\
Set the $\tilde{\mathcal S}_n$ as defined in (\ref{eq:S}).\\
Compute $u_{n+1}$, i.e., $u_{n+1} := \displaystyle \argmin_{u \in \X} \Bigg[\tilde{\mathcal S}_n \left(\taylor{u_{n}}{u}\right) + \alpha_n \norm{u - u_0}{\X}^2\Bigg].$\\
Update the regularization parameter $\alpha_{n+1}$.
}
\Output{$u_{M}$.}
\end{algorithm}

\subsection{On the randomized sketching operator}

In what comes next,  we use the following notations to easy the readability of the paper: 
\begin{align}
\label{g_n_dagger}
   F_{n}^{\dagger} &:=\taylor{u_n}{\udag};  
   \end{align}
   \begin{align}
   \label{g_n_nplus1}
      F_{n}^{n+1} &:= \taylor{u_n}{u_{n+1}}.
\end{align}
Note that in the case $\mathcal{P}_n$ is the identity (meaning that we use full batch in the case of sampling), we have $\cSt{
\cdot} = \cS{\cdot}$. 
The sketching operator $\mathcal{P}_n$ plays a key role in our proposed framework. 
For a better illustration, we can consider the case where the computation of $F'(u_{n})$ is very expensive, not feasible numerically or only a part of it can be computed.  In this case, the  operator ${\mathcal P}_n$  can be represented by a diagonal matrix $P_n$ with ones and zeros elements randomly distributed
over its diagonal. We can also consider the sketched version of the Jacobean, in the definition of $\mathcal{P}_n$, i.e, $P_n\, F'(u_{n})$.
Note that for the new matrix $P_n\, F'(u_{n-1})$, we need to compute only the derivatives of the function $F$ component for which we have one in the corresponding position in  $P_n$. Therefore, by controlling the "number of ones" in $P_n$, we control the number of components of $F$ for which we compute the derivative. Therefore, in this case, the sketching is equivalent to consider a random mini batch of the sum instead of considering the whole sum at each iteration.

Another example is where the precise evaluation of the operator $F$ and its Jacobian are unattainable. This situation is particularly evident when, given a specific input $x$, the operator $F(x)$ is defined as the expected value over a random variable $\xi$, for instance $F(x) = \mathbb{E}_{\xi}(\tilde F(x,\xi))$. In such cases, the only viable approach for estimating $F$ and its Jacobian is through sampling.  In this work we will make the following assumption on the sketching operator. The assumption is inspired by the standard mini-batch assumption in the stochastic gradient method~\cite{HRS16}.
 \begin{ass} \label{asm:1}
    The randomized sketching operator satisfies the following properties:
    \begin{itemize}
        \item  $
\mathbb{E}\left[\mathcal{P}_n\left[F_{n}^{\dagger}\right]\right]= F_{n}^{\dagger}.
$
\item there exists $\sigma\ge 0$, such that
$\mathbb{E}\left[ \left\|\mathcal{P}_n\left[F_{n}^{n+1}\right]  - F_{n}^{n+1} \right\|\right] \le \sigma$.
    \item there exists a constant $\delta>0$, such that $\|\mathcal{P}_n\left[F_{n}^{n+1}\right] - Z_n\| \le \delta$.

    \end{itemize}
\end{ass}

\subsubsection{Discussion on Assumption \ref{asm:1}}

We note that for the first part of this assumption, the sketching operator when it is applied to $F_{n}^{\dagger}$ is leading to unbiased stochastic estimates for all objectives. These type of assumption are very common in the machine learning community where the sketching operator is often related to the use of subsampling over small fractions of the data~\cite{HRS16}. The second part of the assumption, is related to the fact that the variance of the sketching operator $\mathcal{P}_n$ is bounded which is also a very standard assumption similar when we consider the particular case where the residual function is scalar and sketching operator operates via subsampling~\cite{HRS16}. The third part of the assumption indicates that the gap between the $\mathcal{P}_n\left[F_{n}^{n+1}\right]$ and $Z_n$ has to be uniformly bounded. In practice, $Z_n$ is often selected as the average of observations over $F_{n}^{n+1}$.

\section{Convergence Analysis}
\label{sec:theory_det}
{In this section we present our main theoretical contribution which is a convergence analysis of our SIRGNM method. This will require a preliminary error decomposition, based on different error terms. We will then discuss our main result, which provides a bound on the expected residual, independent of $\alpha_n$.

\subsection{Preliminary error decomposition}

In this section, we firstly begin by introducing some notation we will throughout. Specifically, we introduce the effective noise level as follows:
\[
\err\left(g\right):= \cT{g} - \left(\cS{g} - \cS{\ydag}\right), \qquad g \in \Y,
\]
and its stochastic version 
\begin{align*}
\tildeerr\left(g\right)&:= \cT{g} - \left(\cS{g} - \cS{\ydag}\right)+  \left(\cS{g} - \cSt{g}\right)\\
&=
\err\left(g\right) +  \left(\cS{g} - \cSt{g}\right), \qquad g \in \Y.
\end{align*}

Let also 
\begin{align}
{\partial \tildeerr} &:= \tildeerr\left(F_{n}^{n+1}\right) - \tildeerr\left(F_{n}^{\dagger}\right), \label{eq:errntild}\\
{\partial \err} &:= {\err}\left(F_{n}^{n+1}\right) - {\err}\left(F_{n}^{\dagger}\right), \label{eq:errn} \\
{\partial\mathcal{S}}^{n+1}_n &:= \cS{F_{n}^{n+1}} - \cSt{F_{n}^{n+1}} \label{eq:errs},\\
{\partial\mathcal{S}}^\dagger_n & :=  \cS{F_{n}^{\dagger}} - \cSt{F_{n}^{\dagger}}  \label{eq:errsn}.
\end{align}
 Then, we have  
 \begin{align*}
{\partial \tildeerr} &= {\partial \err} + {\partial\mathcal{S}}^{n+1}_n +{\partial\mathcal{S}}^\dagger_n,  
\end{align*}
The first term in $ {\partial \tildeerr} $ is due to the noise in the data and in case of noise-free observations it vanishes. The rest is due to sampling.

Now, using Assumption~\ref{asm:1}, we introduce the following lemma, which  follows naturally.
\begin{lemma}\label{lem:unbiasederror}
 Let Assumption~\ref{asm:1} hold. Then the  randomized error $\tildeerr$ evaluated at $F_{n}^{\dagger}$ is an unbiased estimator of ${\err}$ at $F_{n}^{\dagger}$, i.e.,
 $$
\mathbb{E}[\tildeerr(F_{n}^{\dagger})]= \err(F_{n}^{\dagger}).$$
\end{lemma}
Then next auxiliary results is usefull for our convergence analysis.
\begin{lemma}\label{lem:min_analysis}
Under our working assumptions,  
one has
\begin{equation*}
\alpha_n \left[ \norm{u_{n+1} - {u}_0}{\X}^2 - \norm{\udag - {u}_0}{\X}^2\right]  + \frac12  \norm{F_{n}^{n+1} - \ydag}{\Y}^2 \leq {\partial \tildeerr} + \frac12 \norm{F_{n}^{\dagger} - \ydag}{\Y}^2. 
\end{equation*}
\end{lemma}
\begin{proof}

 Indeed, one has,
\begin{align*}
&\cSt{F_{n}^{\dagger}} - \cSt{F_{n}^{n+1} }\\
=& \left(\cSt{F_{n}^{\dagger}} - \cS{F_{n}^{\dagger}} \right)+ \left(\cS{F_{n}^{\dagger}} - \cS{\ydag}\right)\\
& +\left(\cS{\ydag} - \cS{F_{n}^{n+1}} \right) + \left( \cS{F_{n}^{n+1} }- \cSt{F_{n}^{n+1}}\right) \\
=& \left({\err}\left(F_{n}^{\dagger}\right) - \tildeerr\left(F_{n}^{\dagger}\right) \right)+ \left(\cT{F_{n}^{\dagger}}- {\err}\left(F_{n}^{\dagger}\right) \right)\\
& - \left(\cT{F_{n}^{n+1}}- {\err}\left(F_{n}^{n+1}\right) \right) + \left( \tildeerr\left(F_{n}^{n+1}\right) - {\err}\left(F_{n}^{n+1}\right)\right) 
\end{align*}
Hence, 
\begin{equation} \label{lm:eq1}
 \cSt{F_{n}^{\dagger}} - \cSt{F_{n}^{n+1} } 
  =  \frac12 \norm{F_{n}^{\dagger} - \ydag}{\Y}^2  - \frac12 \norm{F_{n}^{n+1} - \ydag}{\Y}^2 + \tildeerr\left(F_{n}^{n+1}\right)  - \tildeerr\left(F_{n}^{\dagger}\right).
\end{equation}
On the other hand, the minimality condition of \eqref{eq:var_sgn} implies that 
 \begin{equation*}
 \cSt{F_{n}^{n+1}} + \alpha_n \norm{u_{n+1} - u_0}{\X}^2  \leq \alpha_n  \norm{\udag - u_0}{\X}^2 +\cSt{F_{n}^{\dagger}}. \label{eq_min_analysis_1}
\end{equation*}
Hence, 
\begin{equation*}
\alpha_n \left[ \norm{u_{n+1} - u_0}{\X}^2 - \norm{\udag - u_0}{\X}^2\right] \leq \cSt{F_{n}^{\dagger}} - \cSt{F_{n}^{n+1}}. 
\end{equation*}
Now, by using the equality \eqref{lm:eq1}, one gets
\begin{equation*}
\alpha_n \left[ \norm{u_{n+1} - {u}_0}{\X}^2 - \norm{\udag - {u}_0}{\X}^2\right]  + \frac12  \norm{F_{n}^{n+1} - \ydag}{\Y}^2 \leq {\partial \tildeerr} + \frac12 \norm{F_{n}^{\dagger} - \ydag}{\Y}^2. \label{eq_min_analysis_2}
\end{equation*}
\end{proof}
To proceed further, we need the following variational source condition, which has been first formulated in \cite{hkps07} and has become a standard assumption in the analysis of variational regularization methods. In many situations it turns out that variational source conditions are necessary and sufficient for convergence rates \cite{hw17}. Note that - as typical for source conditions in general - the smoothness of $\udag$ is therein measured relative to the smoothing properties of $F$.

\begin{ass}[Variational source condition]\label{ass:vsc}
There exists a concave index function $\varphi$ (i.e. $\varphi(0) = 0$ and $\varphi \nearrow 
\infty$) such that for all $u \in D(F)$ it holds
\begin{equation}\label{eq:vsc}
\norm{u-\udag}{\X}^2 \leq \norm{u-\hat{u}_0}{\X}^2 - \norm{\udag - \hat{u}_0}{\X}^2 + \varphi \left(\frac12 \norm{\op(u) - \op\left(\udag\right)}{\Y}^2\right).
\end{equation}
\end{ass}
In order to further treat the nonlinearity, we employ also the
following assumption.
\begin{ass}[Tangential cone condition]\label{ass:tcc}
There exists constants $\Ctc \geq 1$ and $\eta > 0$ sufficiently small such that
\begin{align*}
\frac{1}{\Ctc} \norm{\op(v) - \ydag}{\Y}^2 - \eta \norm{\op(u) - \ydag}{\Y}^2 & \leq  \norm{\taylor{u}{v} - \ydag}{\Y}^2\\
& \leq \Ctc \norm{\op(v) - \ydag}{\Y}^2 + \eta \norm{\op(u) - \ydag}{\Y}^2 .
\end{align*}
\end{ass}
\begin{rem}
This tangential cone condition follows from the standard tangential cone condition with some $\Ctc$, see \cite[Lemma 5.2]{HW13}. If $\varphi \geqslant \sqrt{t}$ as $t \to 0$, it can - using the techniques from \cite{FW15} - be replaced by a Lipschitz-type assumption.
\end{rem}
In what comes next, we will use the following notation: 
\begin{align}
d_n &:= \mathbb{E}\left[\norm{u_n - \udag}{\X}^2\right],\\
t_n &:= \frac12\mathbb{E}\left[\norm{\op\left(u_n\right) - \ydag}{\Y}^2 \right].
\end{align}
 Therewith, we have proven the following:
\begin{Lem}[Preliminary error estimate]
\label{lem:free}
Let Assumptions \ref{ass:vsc} and \ref{ass:tcc} hold and assume that $u_n \in D(F)$ is well defined. Then  $u_{n+1}$  satisfies  the  following  error decomposition
\begin{equation}\label{eq:aux2}
	\alpha_n d_{n+1} + \frac{1}{2\Ctc} t_{n+1} \leq \mathbb{E}\left[{\partial \tildeerr}\right] + \alpha_n \Psi \left(2 \Ctc \alpha_n\right) + 2 \eta t_n,
\end{equation}
where we abbreviate $\Psi\left(2 \Ctc \alpha_n \right)$ as
\begin{align}
\Psi\left(2 \Ctc \alpha_n \right)& := \left(-\varphi\right)^* \left(-\frac{1}{2 \Ctc \alpha_n}\right) = \sup_{\tau \geq 0} \left[\varphi\left(\tau\right) - \frac{\tau}{2 \Ctc \alpha_n}\right], 
\end{align}
with  $\left(-\varphi\right)^*$  being the Fenchel conjugate of the convex function $-\varphi$.
\end{Lem}
\begin{proof}
Indeed, plugging Assumption \ref{ass:vsc} into Lemma~\ref{lem:min_analysis} with $u = u_{n+1}$ yields
\begin{multline}\label{eq:aux}
\alpha_n \norm{u_{n+1} - \udag}{\X}^2 + \frac12  \norm{F_{n}^{n+1} - \ydag}{\Y}^2 \\\leq \tildeerr + \alpha_n\varphi \left(\frac12 \norm{\op\left(u_{n+1}\right) - \op\left(\udag\right)}{\Y}^2\right) + \frac12 \norm{F_{n}^{\dagger} - \ydag}{\Y}^2.
\end{multline}
   Now, using $v=u_{n+1}$ and $u=u_n$ in the tangential cone condition (from Assumption~\ref{ass:tcc}) gives for the second term on the left-hand side of \eqref{eq:aux} that
\[
\frac{1}{\Ctc} \norm{\op\left(u_{n+1}\right) - \ydag}{\Y}^2 - \eta \norm{\op\left(u_n\right) - \ydag}{\Y}^2 \leq \norm{F_{n}^{n+1} - \ydag}{\Y}^2,
\]
and for the third term on the right-hand side, with (\ref{eq:ip}), that
\begin{align*}
\norm{F_{n}^{\dagger}-\ydag}{\Y}^2 &\leq \Ctc \norm{\op\left(\udag\right) - \ydag}{\Y}^2 + \eta \norm{\op\left(u_n\right) - \ydag}{\Y}^2 = \eta \norm{\op\left(u_n\right) - \ydag}{\Y}^2.
\end{align*}
By inserting the last two inequalities (above) into (\ref{eq:aux}), we obtain the recursive error estimate
\begin{multline}
\alpha_n \norm{u_{n+1} - \udag}{\X}^2 + \frac{1}{2\Ctc}  \norm{\op\left(u_{n+1}\right) - \ydag}{\Y}^2 \\\leq  {\partial \tildeerr} + \alpha_n\varphi \left(\frac12 \norm{\op\left(u_{n+1}\right) - \op\left(\udag\right)}{\Y}^2\right) + \eta  \norm{\op\left(u_n\right) - \ydag}{\Y}^2. \label{eq_min_analysis_3_}
\end{multline}
Now by taking the expectation on both sides of the last inequality (and then apply the Jensen's inequality), one gets
\begin{multline*}
\alpha_n \mathbb{E}\left[\norm{u_{n+1}- \udag}{\X}^2\right]  + \frac{1}{2\Ctc} \mathbb{E}\left[ \norm{\op\left(u_{n+1}\right) - \ydag}{\Y}^2\right] \\\leq \mathbb{E}\left[ {\partial \tildeerr}\right] + \alpha_n\varphi \left(\frac12 \mathbb{E}\left[ \norm{\op\left(u_{n+1}\right) - \op\left(\udag\right)}{\Y}^2\right]\right) + \eta \mathbb{E} \left[\norm{\op\left(u_n\right) - \ydag}{\Y}^2\right]. \label{eq_min_analysis_3}
\end{multline*}
Equivalently,
\begin{equation}
\alpha_n d_{n+1}  + \frac{1}{\Ctc} t_{n+1} \leq \mathbb{E}\left[ {\partial \err}\right] + \alpha_n\varphi \left(t_{n+1}\right) + 2 \eta t_n. \label{eq_min_analysis_4}
\end{equation}
On the other hand, by using the properties of the Fenchel conjugate of $-\varphi$,
\begin{align*}
\Psi\left(2 \Ctc \alpha_n \right)& =  \left(-\varphi\right)^* \left(-\frac{1}{2 \Ctc \alpha_n}\right) = \sup_{\tau \geq 0} \left[\varphi\left(\tau\right) - \frac{\tau}{2 \Ctc \alpha_n}\right] \ge  \varphi \left(t_{n+1}\right) - \frac{t_{n+1}}{2\Ctc \alpha_n}.
\end{align*}
Thus,
\begin{equation}
    \alpha_n\varphi \left(t_{n+1} \right) \le \alpha_n \Psi\left(2 \Ctc \alpha_n \right) + \frac{t_{n+1}}{2\Ctc}.
    \label{eq_min_analysis_5}
\end{equation}

Hence, by combining inequalities \eqref{eq_min_analysis_4} and \eqref{eq_min_analysis_5}, the proof is concluded.
\end{proof}

\begin{lemma}
\label{lem:free2}
Under Assumption \ref{asm:1}, we have that
   $$ \mathbb{E}[{\partial \tildeerr}] \leq \mathbb{E}[\partial \err] + \frac{\sigma^2}{2} + \delta \sigma.$$ 
\end{lemma}
\begin{proof}
  We have 
   \begin{align*}
{\partial \tildeerr} &= {\partial \err} + {\partial\mathcal{S}}^{n+1}_n +{\partial\mathcal{S}}^\dagger_n,  
\end{align*}
Since $F_{n}^{\dagger}$ does not depend on $u_{n+1}$, thus by Lemma \ref{lem:unbiasederror}, one has 
$$
\mathbb{E}[{\partial\mathcal{S}}^\dagger_n]=
\mathbb{E}\Big[\cSt{F_{n}^{\dagger}}-\cS{F_{n}^{\dagger}}\Big] = 0,
$$
Note that for the remaining term \eqref{eq:errn}, i.e. ${\partial\mathcal{S}}^{n+1}_n$, we can not use Lemma \ref{lem:unbiasederror} because of the dependence between $F_{n}^{n+1}$, $u_{n+1}$ and the sketching operator $\mathcal{P}_n$. In this case, one has
\begin{align*}
    {\partial\mathcal{S}}^{n+1}_n &= \cS{F_{n}^{n+1}} - \cSt{F_{n}^{n+1}}\\
    &= \frac{1}{2}\|F_{n}^{n+1}\|^2_{\mathcal{Y}_p}- \left<Z_n,  F_{n}^{n+1}\right> - \frac{1}{2}\|\mathcal{P}_n\left[F_{n}^{n+1}\right]\|^2_{\mathcal{Y}_p}+ \left< Z_n,  \mathcal{P}_n\left[F_{n}^{n+1}\right]\right>\\
     &= \frac{1}{2}\|F_{n}^{n+1} - \mathcal{P}_n\left[F_{n}^{n+1}\right] \|^2_{\mathcal{Y}_p} + \left<\mathcal{P}_n\left[F_{n}^{n+1}\right]-Z_n,  F_{n}^{n+1} - \mathcal{P}_n\left[F_{n}^{n+1}\right] \right> \\
     &\le  \frac{1}{2}\|F_{n}^{n+1} - \mathcal{P}_n\left[F_{n}^{n+1}\right] \|^2_{\mathcal{Y}_p} + \left\|\mathcal{P}_n\left[F_{n}^{n+1}\right]-Z_n  \right\|_{\mathcal{Y}_p}   \left\|F_{n}^{n+1} - \mathcal{P}_n\left[F_{n}^{n+1}\right] \right\|_{\mathcal{Y}_p}.
\end{align*}
Hence using Assumption \ref{asm:1}, one gets
$$
\mathbb{E}[{\partial\mathcal{S}}^{n+1}_n] \le 
\frac{\sigma^2}{2} + \delta\sigma.
$$
This concludes the proof.
\end{proof}

For the sake of readability and simplicity of our final results on the convergence rates and without loss of generality, we assume that $\varphi(t) = \sqrt{t}$, therefore 
$\psi(t) = t/2$.

The following theorem provides the  main bound on the expected residual at iteration $n+1$. This bound comprises three terms: the first term accounts for the noise error, the second term can be managed by selecting an appropriate step size rule, and the third term is controlled by choosing a sufficiently small parameter $\eta$. 


\begin{thm}
\label{thm:mainthm}
 Let Assumptions \ref{asm:1}, \ref{ass:vsc} and \ref{ass:tcc} hold. Assume  that $\eta < \frac{1}{4 \Ctc}$ then we have the following convergence rates
\begin{align*}
t_{n+1} 
 &\leq  \left(2\Ctc e_{\max}  +\Ctc \sigma^2 +2\Ctc \delta \sigma\right)\sum_{i = 0}^n (4 \Ctc \eta)^i  + 
 2 \Ctc^2  \sum_{i = 0}^n \left(4 \Ctc \eta\right)^i \alpha_{n-i}^2
 + \left(4 \Ctc \eta\right)^n t_0,
 \end{align*}
    where $e_{\max} = \max_{n\ge 0}\mathbb{E}\left[{\partial \err}\right] $.
\end{thm}
\begin{proof}
\item From Lemma \ref{lem:free} and Lemma \ref{lem:free2}, 
we have
\begin{align*}
t_{n+1} &\leq 2\Ctc \mathbb{E}\left[{\partial \tildeerr}\right] + 2 \Ctc^2 \alpha_n^2 + 4 \Ctc \eta t_n \\
 &\leq 2\Ctc \mathbb{E}\left[{\partial \err}\right] + \Ctc \sigma^2 + 2\Ctc \delta \sigma + 2 \Ctc^2 \alpha_n^2 + 4 \Ctc \eta t_n.
\end{align*}
By unrolling the recurrence we get 
\begin{align*}
t_{n+1} 
 &\leq  \left(2\Ctc e_{\max}  +\Ctc \sigma^2 +2\Ctc \delta \sigma\right)\sum_{i = 0}^n (4 \Ctc \eta)^i  + 
 2 \Ctc^2  \sum_{i = 0}^n \left(4 \Ctc \eta\right)^i \alpha_{n-i}^2
 + \left(4 \Ctc \eta\right)^n t_0.
 \end{align*}
\end{proof}
The following corollarly gives the convergence rates for the constant and the geometrically decreasing regularization.  
\begin{cor}
\label{thm:main}
Let Assumptions \ref{asm:1}, 
\ref{ass:vsc} and 
\ref{ass:tcc} hold. Assume also that we have a constant or linear decreasing step size, i.e. there exists  
$\gamma \leq 1$ such that $\alpha_n = \gamma \alpha_{n-1}$ and that $\eta < \frac{\gamma^2}{4 \Ctc}$ then we have the following convergence rates
\begin{enumerate}
    \item (general case)
    $$ t_{n+1} \leq  \frac{\left(2\Ctc e_{\max}  +\Ctc \sigma^2 +2\Ctc \delta \sigma\right)}{1 - 4 \Ctc \eta} +  \mathcal{O}(\alpha_n^2),$$
    \item (noise free case, and full sampling i.e, $\partial \err =0$) 
    $$ t_{n+1} =    \mathcal{O}(\alpha_n^2) \text { and } d_{n+1} =  \mathcal{O}(\alpha_n). $$
\end{enumerate}
\end{cor}
\begin{proof}

\begin{enumerate}
\item From Theorem  \ref{thm:mainthm}   we have
\begin{align*}
t_{n+1} &\leq  \left(2\Ctc e_{\max}  +\Ctc \sigma^2 +2\Ctc \delta \sigma\right) \sum_{i = 0}^n (4 \Ctc \eta)^i  + 
 2 \Ctc^2  \sum_{i = 0}^n \left(4 \Ctc \eta\right)^i \alpha_{n-i}^2
 \\&+ \left(4 \Ctc \eta\right)^n t_0. \end{align*}
 Therefore by using the fact that $\alpha_n = \gamma \alpha_{n-1}$ we get
 \begin{align*}
 t_{n+1} &\leq  \left(2\Ctc e_{\max}  +\Ctc \sigma^2 +2\Ctc \delta \sigma\right)\sum_{i = 0}^n (4 \Ctc \eta)^i  + 2 \Ctc^2 \alpha_0^2 \gamma^{2n}  \sum_{i = 0}^n \left(\frac{4 \Ctc \eta}{\gamma^2}\right)^i  \\&+ \left(4 \Ctc \eta\right)^n t_0\\
&\leq   \frac{\left(2\Ctc e_{\max}  +\Ctc \sigma^2 +2\Ctc \delta \sigma\right)}{1 - 4 \Ctc \eta} +  \frac{2 \Ctc^2 \alpha_0^2 \gamma^{2n}}{1 - \frac{4 \Ctc \eta}{\gamma^2}} + \left(4 \Ctc \eta\right)^n t_0 \\
&= \frac{\left(2\Ctc e_{\max}  +\Ctc \sigma^2 +2\Ctc \delta \sigma\right)}{1 - 4 \Ctc \eta} +  \mathcal{O}(\alpha_n^2).
\end{align*}

\item  From Lemma \ref{lem:free} we have
\begin{equation*}
	 d_{n+1} \leq  2 \Ctc^2 \alpha_n + \frac{4 \Ctc \eta t_n}{\alpha_n} - \frac{t_{n+1}}{\alpha_n \Ctc} = \mathcal{O}(\alpha_n).
\end{equation*}
\end{enumerate}
\end{proof}

The previous corollary establishes that the residuals exhibit convergence toward a neighborhood of zero, and this convergence occurs at a rate proportional to $\alpha_n^2$. In other words, as the sequence progresses, the residuals progressively approach a neighborhood of zero with quadratic rate.
In the noise free case, the theorem asserts the dual convergence of both the residual and iterate sequences toward zero. Specifically, the residual sequence converges with a rate proportional to $\alpha_n^2$, indicating a quadratic decrease, while the iterate sequence exhibits a convergence rate of $\alpha_n$, implying a linear decrease.

\section{Numerical Experiments}
\label{sec:num}

In this section we provide numerical experiments demonstrating, and highlighting the benefit, of our stochastic IRGNM. We apply this in the context of both noise-free and noisy data, where we test our algorithm on a random PDE example, taking motivation from geophysical sciences. We will compare our methodology to previously discussed methods which demonstrate an improvement on accuracy. 

\subsection{Overview of algorithms}
\textcolor{black}{We briefly provide an overview of our algorithms before discussing our numerical simulations. Throughout the experiments we use Algorithms \ref{alg:cIRGNM}, \ref{alg:SIRGNM} and \ref{alg:SdIRGNM}. Algorithm~\ref{alg:cIRGNM} consists of the version of the vanilla iterative scheme, referred to as the IRGNM. For this algorithm we  consider the ''classical" version which contains no modified noise, where $Z_n=Y$ as from \eqref{eq:model_dyn}, and the dynamic version with sequential observations from \eqref{eq:avg_noise}.
 This takes motivation from the work of Chada et al.~\cite{CIL22}, which takes sequential observations. We  then consider a stochastic version of this which we refer it as the sIRGNM in Algorithm~\ref{alg:SIRGNM}. 
Finally our last algorithm is Algorithm~\ref{alg:SdIRGNM} which is the stochastic version of the dynamic algorithm, referred to as SdIRGNM.} 

\subsection{Darcy flow model}

We consider the inversion of the Darcy flow model~\cite{MH97}, which models groundwater flow
in a porous medium. For this we consider the problem of solving an elliptic partial differential equation (PDE), 
where  the forward problem is to solve the PDE 
\begin{equation}
\label{eq:darcy1}
-\nabla\cdot({\kappa}\nabla p) = f,   \ x \ \in \ \Omega,
\end{equation}
with mixed boundary conditions
\begin{equation}
p(x_1,0) = 100, \ \ \ \  \frac{\partial p}{\partial x_1}(6,x_2) = 0, \ \ -\kappa\frac{\partial p}{\partial x_1}(0,x_2) = 500, \ \ \  \frac{\partial p}{\partial x_2} (x_1,6) = 0,
\end{equation}
and the source term $f$ defined as
\[
f(x_1,x_2)=
\begin{cases}
0, & \textrm{if} \ \  0 \leq x_2 \leq 4, \\
137, & \textrm{if} \ \  4 \leq x_2 \leq 5, \\
274,  & \textrm{if} \ \  5 \leq x_2 \leq 6,
\end{cases}
\]
where $\kappa \in L^{\infty}(\Omega)$ denotes the permeability, $f \in L^{\infty}(\Omega)$ is a source term and $p \in H^1(\Omega)$ is the pressure, representing the forward solution. For our problem the domain is specified as $\Omega = [0,6]^2$.
The inverse problem associated with \eqref{eq:darcy1} is to determine the permeability $\kappa \in L^{\infty}(\Omega)$ from pointwise measurements of the PDE solution  {$p \in H^1(\Omega)$}.
The specific form of the PDE model \eqref{eq:darcy1}, (in relation to the boundary conditions and $f$), is  tested by Hanke \cite{MH97} whom motivated a regularized Levenberg-Marquardt algorithm. More information on this PDE model can be found in \cite{CN86}.  
To treat a more general setting, we consider the following weighted $L^2$ space
\begin{align}
\label{eq:darcy00}
\mathcal{H}\equiv \{u\in L^{2}(\Omega)\big\vert\quad \vert\vert  \mathcal{C}^{-1/2} u \vert \vert_{L{^2}(\Omega) } \leq \infty\},
\end{align}
such that our covariance operator $\mathcal{C}$ induced by a correlation function given by 
\begin{align}\label{eq:eit100}
\mathcal{C}(u)[x]=\int_{\Omega}u(x')c(x,x')dxdx'.
\end{align}
Specifically we choose a Mat\'{e}rn covariance function which takes the form
$$
\mathcal{C}(x,x') = \frac{\Gamma(\nu)}{2^{1-\nu}} K_{\nu}\bigg(\frac{|x-x'}{\ell}\bigg)\bigg(\frac{|x-x'}{\ell}\bigg)^{\nu},
$$
where $(\nu, \ell) \in \bbR^+ \times \bbR$  are
associated hyperparameters which define the smoothness and lengthscale, and $K_{\nu}(\cdot)$ represents a Bessel function of the second kind, 
where use the same notation of $\mathcal{C}$ as the covariance operator, and its corresponding function. We use this function, as it is common when aiming to simulate Gaussian random fields.

We compare the performances of our approaches against their corresponding  deterministic versions.  
The equivalent updates of (\ref{eq:gn}) of the deterministic algorithm cIRGNM \cite{CIL22}
 that take into consideration weighting of $\mathcal{C}$ 
 is given by 
\begin{align}
\label{eq:darcy01}
{u}_{n+1}={u}_{n} + \left( {\op'[{u}_{n}]^*\op'[u_{n}]+\alpha_{n}\mathcal{C}^{-1}}\right)^{-1}{\Big(\op'[{u}_{n}]^*(y-\op({u}_{n}))+\alpha_{n}\mathcal{C}^{-1}({u}_0-{u}_{n})\Big)},
\end{align}
and similarly for the dIRGNM \cite{CIL22}.
For computational efficiency we  use the Woodbury lemma to compute the updates, for instance   \eqref{eq:darcy01} is equivalent to 
\begin{align*}\label{eq:darcy02}
{u}_{n+1}&=u_{0}+{\mathcal{C}\op'(u_{n})^*}{(\op'(u_{n})\mathcal{C}\op'(u_{n})^*+\alpha_{n}I)}^{-1}\Big(y-\op(u_{n})-\op'(u_{n})(u_0-u_{n})\Big). 
\end{align*}
We  consider two experiments with two different groundtruths to recover, the first being a smooth function defined as,
$$
u^{\dagger}(x,y)=\exp\Big[-100\Big( (x-0.3)^2+(y-0.7)^2\big))\Big]+\frac{1}{2}\exp\Big[-100\big( (x-0.7)^2+(y-0.35)^2\big)\Big],
$$
over the domain $\Omega=[0,6]^2$. The second ground truth is a piecewise constant function based on a level set thresholding of a Gaussian random field. 
\newpage
\subsubsection{Parameter choices}
We choose  $64$ pointwise measurements within our domain $\Omega$, and run our experiment for $n=10^4$ iterations. For the regularization parameter, we set {\color{black}$\alpha_0=0.5$}. We modify the choice of the regularization parameter, which is explained in the next subsection, and within the figures. The variance of the data is chosen as $\delta=0.1$. For the Mat\'{e}rn covariance function,  we specify $(\nu,\ell)= (2.4,15)$. 
 $u_0=1$ is chosen to be the vector of ones. 

 For our stochastic algorithm variants, we choose $\mathcal{P}_n$ to be the "minibatch" operator, in the sense that it selects only a part of the vector. Within our experiment, at each iteration, for our stochastic variants (SIRGNM and SdIRGNM) we keep $32$ components selected randomly at uniform. In the numerical experiments, as a stopping criteria, we use the of relative error which is defined as
 $$
\mathrm{err}_{\mathrm{relative}} = \frac{\|u-u^{\dagger}\|_{L^2(\Omega)}}{\| u^{\dagger}\|_{L^2(\Omega)}},
 $$
 which considers the percentage of difference with respect to the $L^2(\Omega)$ norm. \textcolor{black}{For the initial guess of the discontinuous truth we set  $u_0=1$, and $u_0=0.1$ for the smooth ground truth. As discussed in \cite{CIL22}, this choice is not important which was tested for different choices previously.}

\subsection{Summary from figures}
Our first set of plots relate to comparing SIRGNM and SdIRGNM to the deterministic  algorithms (IRGNM and dIRGNM), where we consider the relative $L_2$ error. Figures~\ref{fig:plot2} presents a comparison of different plots where for the noisy-data, we consider alternative choices for the adaptive regularization parameter. As we can observe in both instances the stochastic version performs better than its deterministic counterpart, where we see an improvement in accuracy. The fluctuations that arise, are induced from the uncertainty.  We also provide 95\% confidence intervals, where we notice a low variance where the upper level usually is as good as its deterministic method.

The results obtained both both ground truth  are similar and the stochastic algorithms seem to perform the best. We further show the effect if we consider a constant $\alpha_n$, in comparison with noise-free data cases in Figure \ref{fig:plot3} which further highlights the benefit of the stochastic method, and the lack of significant effect on the constant choice for the regularization parameter. Finally, to get an understanding on the performance of these methods, we compare the final reconstruction of some of these algorithms where we plot the SIRGNM in Figure \ref{fig:11} - \ref{fig:22} and in all cases we see a similar performance, where the noise-free algorithms seem to perform the best. We not the early termination of the noise-free data is based on a threshold, which is attained in both experiments. We also terminate the noise-free simulations early as we aim to demonstrate the convergence speed difference, even though we are primarily interested in noisy observations.


\begin{figure}[h!]
    \centering
      \includegraphics[width=0.48\linewidth]{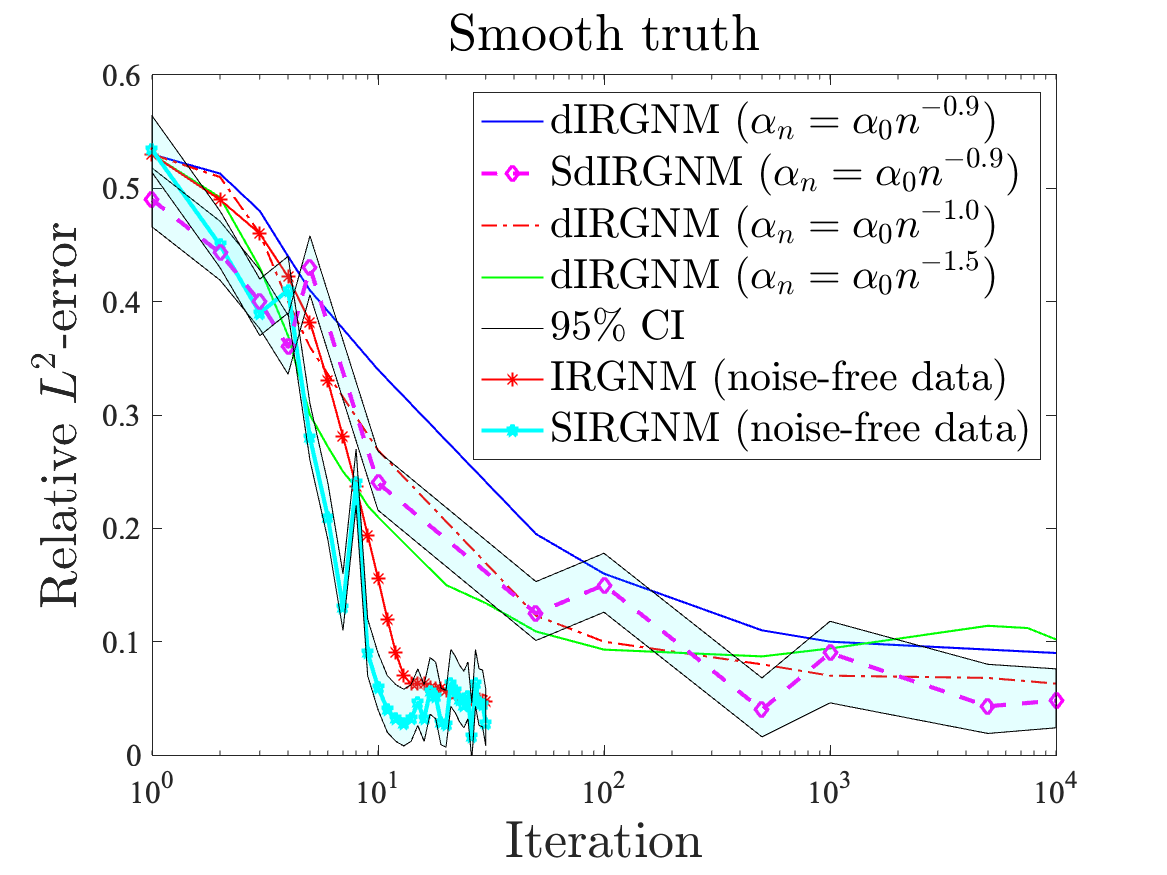}
    \includegraphics[width=0.48\linewidth]{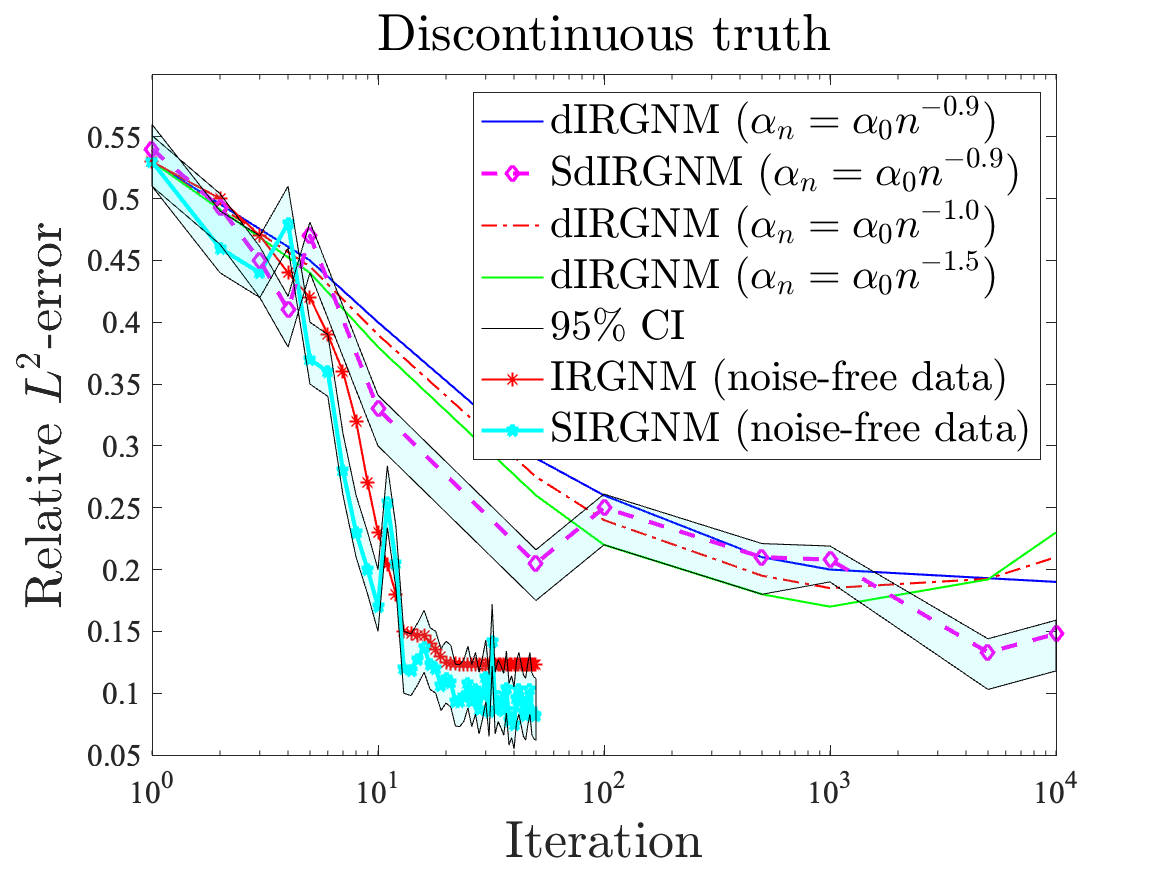}
    \caption{Stochastic IRGNM and dIRGNM for a discontinuous truth $u^{\dagger}$.}
    \label{fig:plot2}
\end{figure}

\begin{figure}[h!]
    \centering
    \includegraphics[width=0.48\linewidth]{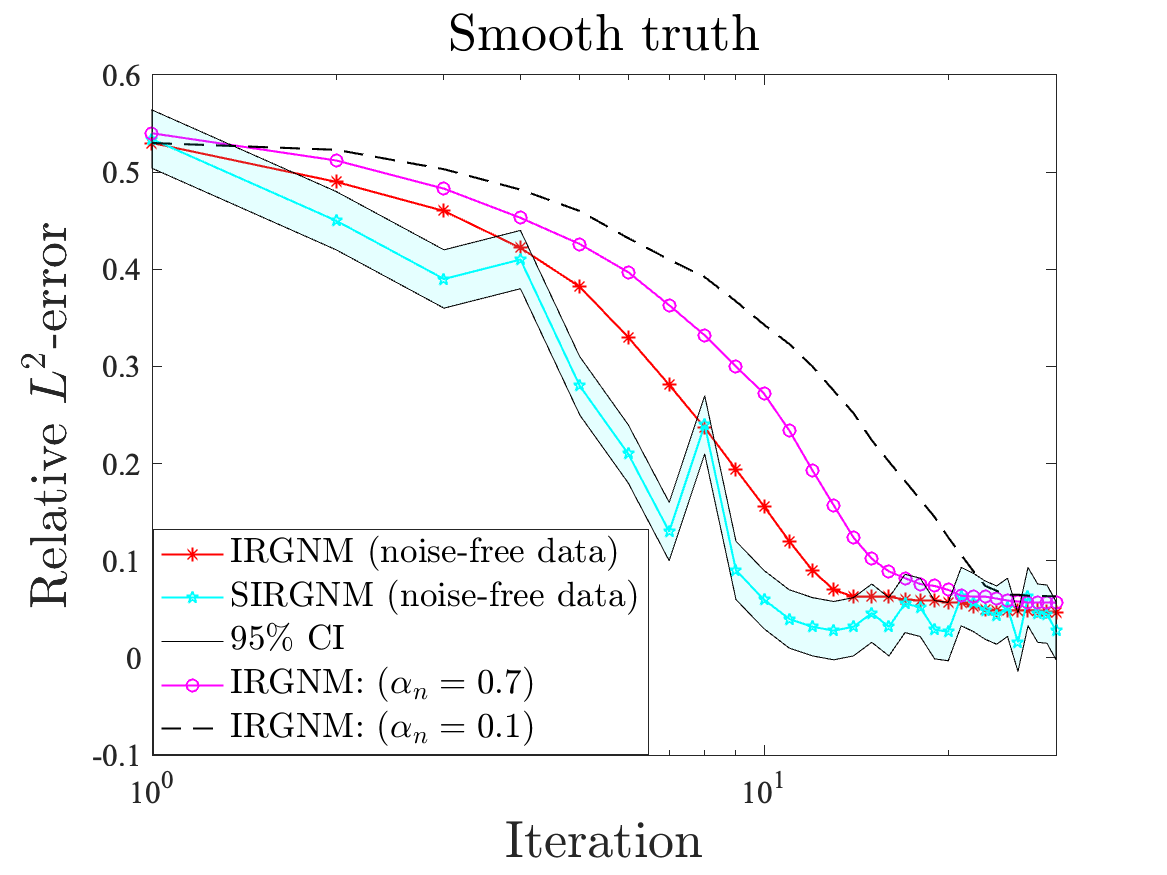}
    \includegraphics[width=0.48\linewidth]{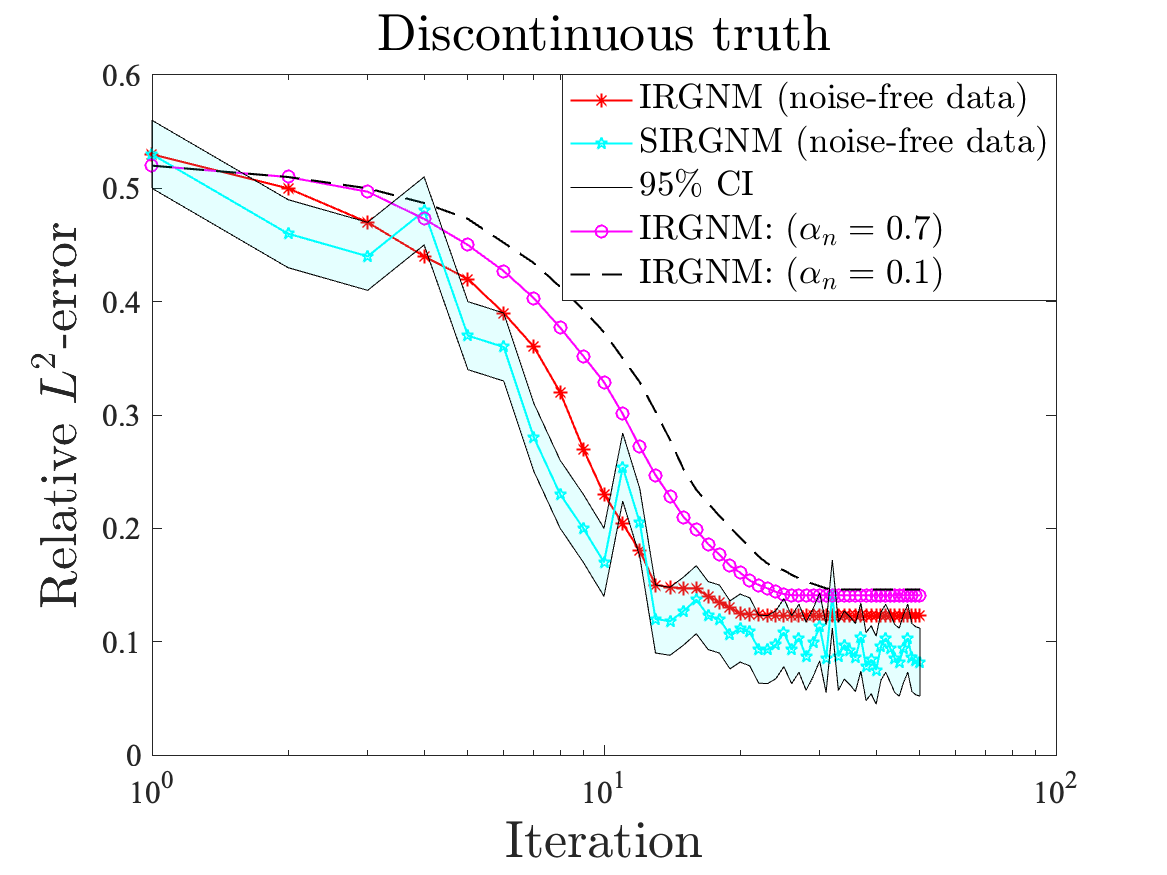}
    \caption{
    IRGNM and SIRGNM
    with constant regularization term.
    }
    \label{fig:plot3}
\end{figure}

\begin{figure}[h!]
\centering
\includegraphics[scale=0.37,trim=0 0 0 0]{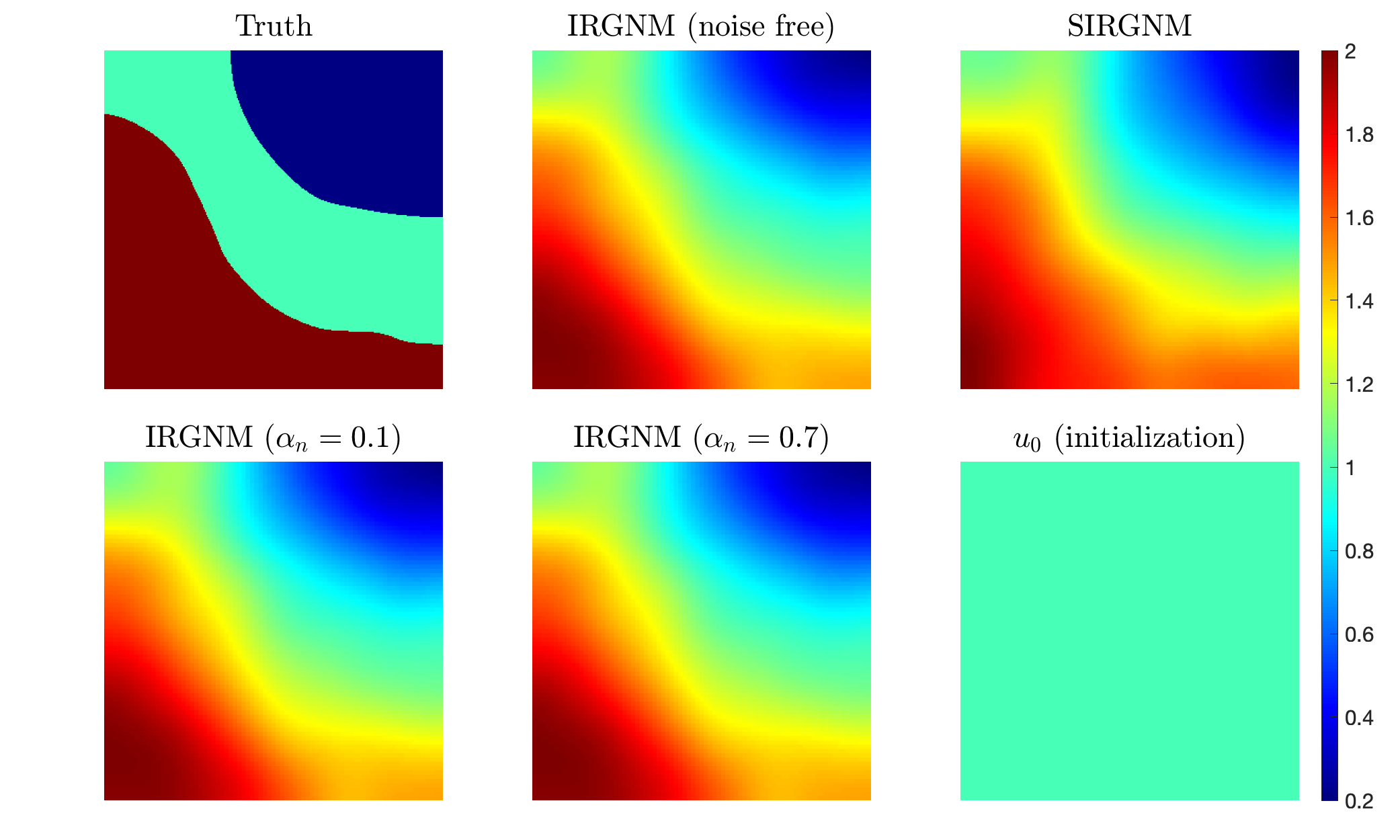}
\caption{Final reconstruction plots of various methods on the discontinuous groundtruth.}
    \label{fig:11}
\end{figure}


\begin{figure}[h!]
\centering
\includegraphics[scale=0.37,trim=0 0 0 0]{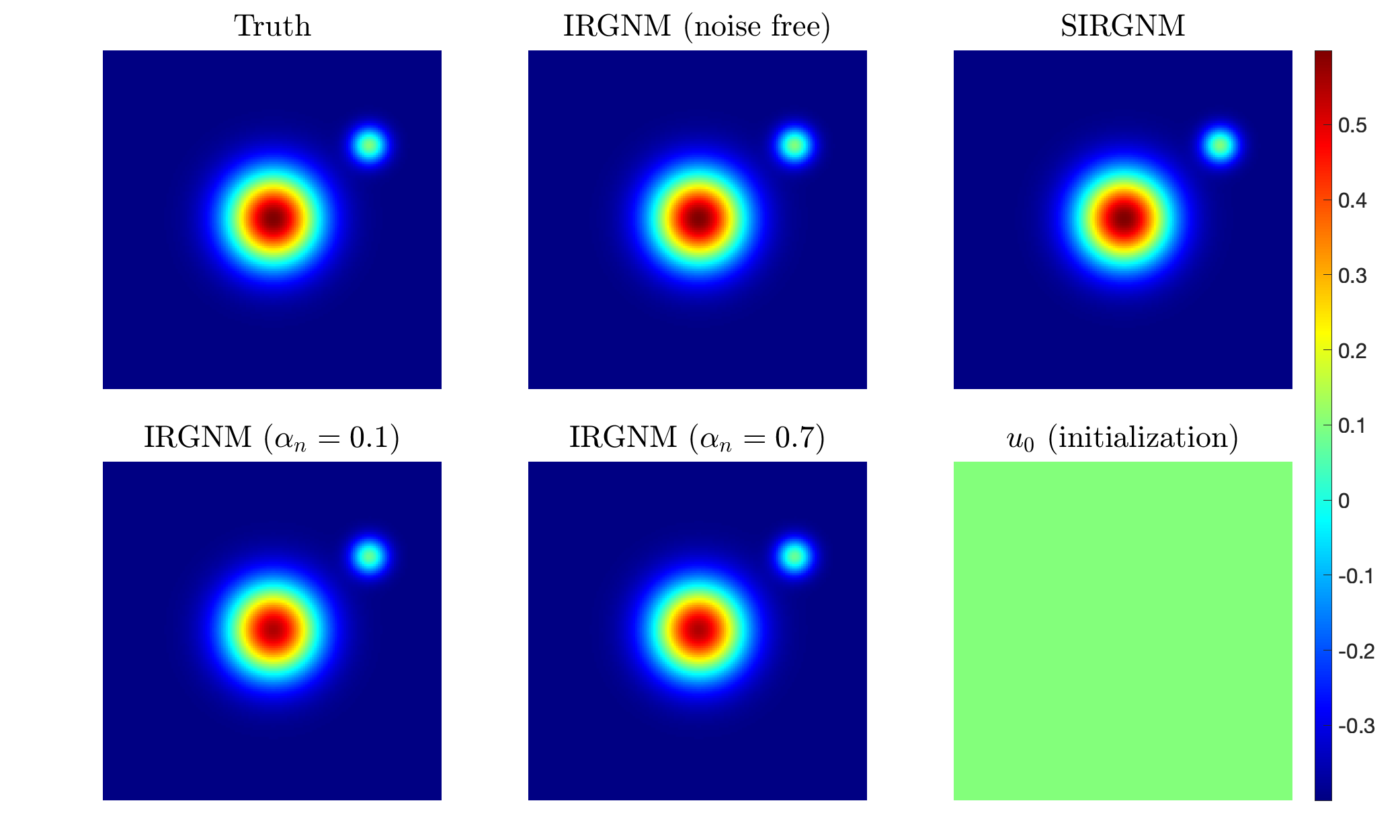}
\caption{Final reconstruction plots of various methods on the smooth groundtruth.}
    \label{fig:22}
\end{figure}

\subsection{Effect of mini-batching on the convergence speed}
Another question, in the context of stochastic optimization algorithms, is the effect of mini-batching. In particular we repeated the experiments above with the effect of mini-batching. Table \ref{table:1} provides results of this, where we change the size of mini-batches corresponding to the data set, and monitor the final reconstruction error at the end of each experiment. Our results here are for both the smooth and non-smooth truths. The table is considered for the adaptive choice of $\alpha_n = \alpha_0 n^{-0.9}$.
\\
As demonstrated in the table, as we increase the batch-size, the CPU time and iterations required to achieve a relative error of $0.1$ increases, which is expected. However the relative errors, independent of the batch-size, does not differ significantly, which improves slightly as we increase the batch-size. 

\begin{table}[]
\begin{tabular}{|l|c|c|c|c|}
\hline
\multicolumn{1}{|c|}{Truth $u^{\dagger}$}            & Mini-batch size & Relative error (final) & \# of iterations &  CPU time\\ \hline \hline
\multirow{4}{*}{Discontinuous  } & 16 &0.198  &  1925 &  4.5 hours \\ \cline{2-5} 
                                       & 32 & 0.173 & 2269 &  4.4 hours \\ \cline{2-5} 
                                       & 64 & 0.172  &  2741  & 4.5 hours \\ \cline{2-5} 
                                       & 128 & 0.166 &  3034 & 4.6 hours \\ \hline \hline
\multirow{4}{*}{Smooth} & 16 &0.092  &  613 & 4.4 hours \\ \cline{2-5} 
                                       & 32 &  0.083 & 697 & 4.4 hours \\ \cline{2-5} 
                                       & 64 &0.077  &  723 & 4.5 hours \\ \cline{2-5} 
                                       & 128 &  0.076 &  856 & 4.7 hours \\ \hline
\end{tabular}
\caption{The effect of mini-batching and truth type on the convergence speed (measured by the number of iterations and related CPU time) of the proposed method. The stopping criterion is set to achieve $\mathrm{err}_{\mathrm{relative}} \le 0.1$.}
\label{table:1}
\end{table}


\section{Conclusion}
\label{sec:conc}
Stochastic gradient-based algorithms have become very important in the fields of applied mathematics and machine learning, where  exact gradient-evaluations can be cumbersome and costly. This argument extends to this work where we develop a stochastic version of the inverse-problem methodology which is the iteratively regularizaed Gauss-Newton method. In particular we also extend this to recent work of Chada et al. \cite{CIL22} which proposed a dynamics version which allows sequential observations.  
We provided the error analysis related to our new approach. 
We then provided numerical experiments on a nonlinear a 2D elliptic PDE example, which demonstrate the improved performance of the SIRGNM method, and its dynamic variant. This is compared to various methods with varying choices of regularization parameters.

In terms of further directions to consider for future work, we list a number of these below.
\begin{itemize}
\item[(i)] This work considered a version of the Gauss-Newton method. It would be of interest to consider to the extension to the regularized Levenberg-Marquardt method, which follows very similarly with a key difference being the penalty term. A more thorough investigation is needed in order to develop sound theory.
\item[(ii)] Another potential direction would be to consider a further analysis for infinite-dimensional optimization such as gradient flow systems.
\item[(iii)]  Using these ideas for sequential observations, and uncertainty could be potentially applied, and used, in understanding the ensemble Kalman filter algorithm for inverse problems \cite{CSW19,CCS21,CIRS18,CT21,WCS22}. Despite the algorithm being derivative-free, it has many connections with both the Gauss--Newton and Levenberg-Marquardt methods \cite{BMB94,BDK22,GLN07}.
\item[(iv)] \textcolor{black}{Another potential work is the setting where one has a corrupt operator, or biased operator with respect to the forward operator $F(\cdot)$. Such a consideration is common in the stochastic optimization which is covered in various pieces of literature \cite{AS20,MKS23}. }
\item[(v)] \textcolor{black}{Finally one could aim to understand various stopping criterion's in the context of stochastic iterative regularization. This has been achieved in the non-dynamic version, however this is still very much in the dynamic setting. One starting point may be the use of the Lepskii principle \cite{LS97,PM06}.}
\end{itemize}

\section*{Acknowledgments}
NKC is supported by an EPSRC-UKRI AI for Net Zero Grant: “Enabling CO2
Capture And Storage Projects Using AI”, (Grant EP/Y006143/1).
NKC is also supported by a City University of Hong Kong startup grant.}

\newpage

\appendix

\section{Algorithms}
\label{app:algo}



\begin{algorithm}[h!]
\caption{Stochastic iterated regularized Gauss-Newton method (\texttt{SIRGNM})}
\label{alg:SIRGNM:1}
\SetKwInOut{Input}{inputs}
 \SetKwInOut{Output}{output}

     \Input{$u_0 \in \X$, $\alpha_0>0$ (an initial value for the regularisation), and $M$ (a maximum number of iterations).}


\For{$n=0,\ldots,M-1$}{
Select  $\mathcal{P}_n$ a stochastic sketching operator.\\
Select  $\alpha_n \in  \mathbb{R}$.\\
Select  $Z_n = Y \in \Y_p$.\\
Set the $\tilde{\mathcal S}_n$ as defined in (\ref{eq:S}).\\
Compute $u_{n+1}$, i.e., $u_{n+1} := \displaystyle \argmin_{u \in \X} \Bigg[\tilde{\mathcal S}_n \left(\taylor{u_{n}}{u}\right) + \alpha_n \norm{u - u_0}{\X}^2\Bigg].$\\
}
    \Output{${u}_{M}$.}
\end{algorithm}

\begin{algorithm}[h!]
\caption{Stochastic dynamic iterated regularized Gauss-Newton method (\texttt{SdIRGNM})}\label{alg:SdIRGNM}
\SetKwInOut{Input}{inputs}
 \SetKwInOut{Output}{outputs}

     \Input{$u_0 \in \X$, $\alpha_0>0$ (an initial value for the regularisation), and $M$ (a maximum number of iterations).}


\For{$n=0,\ldots,M-1$}{
Select  $\mathcal{P}_n$ a stochastic sketching operator.\\
Select  $\alpha_n \in  \mathbb{R}$.\\
Set  $$Z_n =\frac{1}{n}\sum_{i=1}^nY_i\in \Y_p,$$ i.e. the averaged observations.\\
Set the $\tilde{\mathcal S}_n$ as defined in (\ref{eq:S}).\\
Compute $u_{n+1}$, i.e., $u_{n+1} := \displaystyle \argmin_{u \in \X} \Bigg[\tilde{\mathcal S}_n \left(\taylor{u_{n}}{u}\right) + \alpha_n \norm{u - u_0}{\X}^2\Bigg].$\\
}
\Output{${u}_{M}$.}
\end{algorithm}





\end{document}